# COALESCENCE IN A RANDOM BACKGROUND

BY N. H. BARTON[1], A. M. ETHERIDGE[1,2] AND A. K. STURM[3]

*University of Edinburgh, University of Oxford and Technical University Berlin*

We consider a single genetic locus which carries two alleles, labelled $P$ and $Q$. This locus experiences selection and mutation. It is linked to a second neutral locus with recombination rate $r$. If $r = 0$, this reduces to the study of a single selected locus. Assuming a Moran model for the population dynamics, we pass to a diffusion approximation and, assuming that the allele frequencies at the selected locus have reached stationarity, establish the joint generating function for the genealogy of a sample from the population and the frequency of the $P$ allele. In essence this is the joint generating function for a coalescent and the random background in which it evolves. We use this to characterize, for the diffusion approximation, the probability of identity in state at the neutral locus of a sample of two individuals (whose type at the selected locus is known) as solutions to a system of ordinary differential equations. The only subtlety is to find the boundary conditions for this system. Finally, numerical examples are presented that illustrate the accuracy and predictions of the diffusion approximation. In particular, a comparison is made between this approach and one in which the frequencies at the selected locus are estimated by their value in the absence of fluctuations and a classical structured coalescent model is used.

**1. Introduction.** The coalescent process was introduced by Kingman (1982) as a simple and elegant description of the genealogical relationships amongst a set of neutral genes. Although the two theories have developed largely independently, the coalescent is closely related to the classical concept of identity by descent introduced independently by Cotterman and Malécot [see Nagylaki (1989) for a survey]. The original coalescent applies

Received October 2002; revised April 2003.
[1]Supported in part by EPSRC GR/L67899 and BBSRC MM109726.
[2]Supported by an EPSRC Advanced Fellowship.
[3]Supported by The Rhodes Trust and DFG Priority Programme 1033.

to the case of a single panmictic population of constant size, but it extends naturally to describe populations that vary with time or to structured populations in which genes may be found in different places or embedded in different genetic backgrounds. By considering the ancestral selection graph, various forms of selection can also be incorporated [Krone and Neuhauser (1997) and Donnelly and Kurtz (1999)]. However, as the genetic sophistication increases, not only are analytic results unattainable but also the approach becomes increasingly computationally intensive. Moreover, the powerful results of Donnelly and Kurtz rest upon exchangeability of the sample. This means that they lend themselves to describing the genealogy of a random sample from the population. In the problem that we are concerned with here, the sample is *not* random.

The particular problem that we are concerned with is the following. Suppose that selection acts on a single locus which carries two alleles labelled $P$ and $Q$. There is also a strictly positive mutation rate between these two alleles so that neither becomes fixed in the population. The selected locus is linked to a second neutral locus with recombination rate $r$. One can then ask about the genealogy of a sample from the neutral locus. If we *know* the type of each individual in the sample at the selected site, then we have a sample from known locations in a *structured* population, with two demes (determined by the $P$ or $Q$ allele) in which the population size fluctuates *randomly.* Recombination and mutation from $P$ to $Q$ both contribute to migration between the demes, while mutation and selection determine the population sizes. Setting the recombination rate $r = 0$ we recover the case of a single selected locus.

For certain forms of selection (directional and balancing) one can address this problem using the ancestral selection graphs of Krone and Neuhauser (1997). Such graphs trace the lineages of "potential ancestors" of a sample from the selected locus. As one traces backwards in time lineages can "branch" as well as "coalesce." On reaching the most recent common ancestor of all such potential lineages, one then traces back through the graph culling those that in fact did not contribute to the sample. However, this method restricts the form of the selection and, as observed in Przeworski, Charlesworth and Wall (1999), is also computationally demanding, especially if selection is strong, because of the proliferation of potential lineages.

In fact the most common approach to our problem is to assume that fluctuations are sufficiently small that we can approximate the allele frequencies at the selected locus by their value in the absence of fluctuations and then model the genealogy at the neutral locus by the structured coalescent with constant deme sizes [Kaplan, Darden and Hudson (1988), Notohara (1990), Herbots (1997) and Nordborg (1997)]. One might expect this procedure to give good approximations to quantities such as the mean time to the most recent common ancestor of the sample, which is very robust, but not to the variance of the same quantity. In fact, in Section 7 we illustrate that for some parameter values even the mean time to coalescence for a sample of size two can be ill-approximated by this procedure.

Here we adopt an alternative approach, more akin to the classical one. We *retain* the stochastic fluctuations in the two deme sizes. Although we first establish a coalescent-like diffusion approximation for the genealogy of the sample, we actually express the quantities that we are interested in as solutions to a system of ordinary differential equations and our numerical examples will all be obtained by numerically solving this system. In Lemma 3.1 we identify the appropriate "coalescent." There are no surprises here, the model being entirely analogous to a structured coalescent with two demes, but now the jump rates in the coalescent are governed by the Wright–Fisher diffusion that determines the population sizes in each deme. What *is* surprising is the range of parameter values for which the coalescent approximation is valid. One might expect that the mutation rates between alleles $P$ and $Q$ need to be large enough that the allele frequencies stay away from the margins where the diffusion approximation should break down. As we see in Theorem 5.1, in fact we have convergence to the diffusion approximation provided that the mutation rates between the selected alleles are nonzero.

In order to exploit the diffusion approximation, we use it to write down a system of ordinary differential equations for the probability of identity in state for a sample of size two from the neutral site (Theorem 6.1). Predictions of this model are compared not only to those of the Wright–Fisher model that it is approximating, but also to those of the constant deme size model in Section 7. We see in particular that the fluctuations in allele frequency have a significant effect on the probabilities of identity.

This approach was first suggested by Kaplan, Darden and Hudson (1988, 1989) who essentially wrote down the same coalescent approximation, although they did not address mathematical questions of existence of the corresponding process or convergence to the limit. They also wrote down the system of differential equations that determine the total length of the ancestral tree of a sample under this approximation. These are of exactly the same form as the equations for probability of identity that we obtain here. Because of the singular nature of the coefficients and the difficulty in assigning boundary conditions to the system, they develop a novel numerical solution scheme. However, this has not been exploited in the literature. Here we have been able to identify the boundary conditions for the system and as a result our numerical techniques are based on standard software. They are described in detail in the companion paper Barton and Etheridge (2004).

A particular strength of our approach is that it applies to very general forms of selection with essentially no additional computational effort. Numerically it is considerably more efficient than simulations based on the ancestral selection graph and moreover, in contrast to such simulations, the computational efficiency does not decrease as the strength of selection increases. Although we concentrate on the simple setting of a sample of size two embedded in a genetic background with just two possible states, the approach is easily extended to larger samples and more complex genetic backgrounds (albeit at the expense of increased computational complexity); see Remark 6.2. As we see in Section 7, even the simplest context provides considerable scope for investigation of important biological issues. We do not explore it here, but we also see the ordinary differential equations as offering a valuable analytic route to a perturbation analysis of identities in fluctuating backgrounds.

Our primary motivation is the desire to investigate the effects of the selection on the coalescence times for the neutral site. However, we also regard this as an important step in understanding more general versions of the structured coalescent in randomly fluctuating backgrounds. This is crucial to understanding populations with complex spatial or genetic structure where the number of individuals in each background is not sufficiently large that fluctuations can be ignored [see Barton, Depaulis and Etheridge (2002), Barton and Navarro (2002) and references therein].

The rest of the paper is laid out as follows. In Section 2 we describe our model for the (forward in time) evolution of the proportions of the $P$ and $Q$ alleles. Our starting point is a version of the *Wright–Fisher model*. For later comparison with the diffusion approximation, we obtain a system of algebraic equations for the probabilities of identity in state for a sample of size two from such a population. As we explain above Definition 2.1, to obtain the diffusion approximation it is convenient to work with the continuous time counterpart of the Wright–Fisher model, the *Moran model*. For such a model, assuming that the frequencies of the selected alleles, $P$ and $Q$, have reached stationarity, we write down the generator of the process $\{(p^{(1)}(t), n_1^{(1)}(t), n_2^{(1)}(t))\}_{t\geq 0}$ that encodes the backward in time evolution of the selected allele frequencies and the numbers of ancestors of the sample of neutral alleles alive at time $t$ before the present, labeled according to their background ($P$ or $Q$). In Section 3 we rescale the parameters in our model and establish the form of the generator of the corresponding diffusion approximation. The existence of a stochastic process with this generator and convergence of the rescaled processes to this limit are established in Sections 4 and 5, respectively. In Section 6 we write down a system of differential equations for the distribution of coalescence times and hence, for the probability of identity in state in a sample of size two. We establish an iterative solution to the system and indicate the extension to larger samples and more complex genetic backgrounds. In the final section, Section 7, we illustrate, in the case of balancing selection, the extremely good agreement with the probabilities of identity established via the Wright–Fisher model of Section 2. Although in biological applications one is typically concerned with a neutral linked locus, we concentrate on the selected locus ($r = 0$) for easy comparison with alternative approaches. We conclude by comparing the predictions of the model to those obtained by assuming that the frequency of alleles at the selected site is deterministic for different strengths of balancing selection. A full discussion of the biological implications can be found in the companion paper Barton and Etheridge (2004).

**2. The model.** In this section we describe the underlying model. First consider the evolution of frequencies at the selected locus. We write $p$ and $q = 1 - p$ for the proportions of $P$ and $Q$ alleles respectively. Our starting point is a Wright–Fisher model with selection and mutation between types. We assume a diploid population of size $N$. Thus each individual has type $PP$, $PQ$ or $QQ$ and we write $P_{11}$, $P_{12}$ and $P_{22}$ for the corresponding proportions of each type. Then

$$p = P_{11} + \tfrac{1}{2}P_{12},$$
$$q = P_{22} + \tfrac{1}{2}P_{12}.$$

During the reproductive process, each individual has a large (effectively infinite) number of germ cells (cells of the same genotype) that split into gametes (cells containing just one chromosome from each pair). The gametes then fuse at random to form the next generation. We assume that there is selection in favor of certain genotypes. Further there is mutation from type $P$ to $Q$ and vice versa.

Suppose that immediately before the reproductive step, the proportion of type $P$ is $p$. For simplicity we assume *multiplicative* selection. That is, relative fitnesses of $PP : PQ : QQ$ can be expressed in the form $u^2 : uv : v^2$. This means that we can model selection as acting on haploids, so that after selection the proportion of type $P$ will be

$$p^* = \frac{p(1+s)}{1+sp}$$

for some $s$. In the case of directional selection, $s$ is just a constant, but by taking $s$ to be frequency dependent (i.e., a function of $p$) we can approximate more complicated selection acting on the diploid population. For example, balancing selection is modeled by assuming that $s = s_0(p_0 - p)$ for some $0 < p_0 < 1$ and constant $s_0$. If the population size is sufficiently large, this is close to a model of overdominance with relative diploid fitnesses $PP : PQ : QQ$ of $1 - s_0q_0 : 1 : 1 - s_0p_0$, where $q_0 = 1 - p_0$.

We now account for mutation between $P$ and $Q$. Suppose that in each generation a mutation from $P$ to $Q$ has probability $\mu_1$ and from $Q$ to $P$ has probability $\mu_2$. After the mutation step, the proportion of type $P$ is then

$$p^{**} = (1 - \mu_1)p^* + \mu_2(1 - p^*).$$

Finally $2N$ gametes are chosen at random to form the next generation. (These fuse at random into $N$ diploid pairs.) The resulting number of type $P$ chromosomes in the population will then be binomially distributed with $2N$ trials and success probability $p^{**}$.

Now consider the neutral locus. This is on the same chromosome as the selected locus, but we allow for the possibility of *recombination* or *crossover* events. This happens during meiosis (the process of splitting into gametes). We assume that with probability $r$ there is a recombination event between the selected and neutral sites. The result is that the two gametes exchange a portion of chromosome that includes the selected site, but not the neutral site. Consequently, if such an event occurs, for each of the two gametes the portion of the chromosome that includes the neutral locus and that segment including the selected locus come from *different* parental chromosomes. From the point of view of the neutral locus, the two chromosomes swap types at the selected locus. By this mechanism (as well as by mutation) an individual from the sample at the neutral locus can be in a different genetic background from her parent.

So that we can later numerically test the accuracy of the diffusion approximation, we now use the Wright–Fisher model to calculate the probabilities of identity in state at the *neutral locus* of a sample of two genomes whose type at the selected locus is known. (We shall refer to "individuals" in the sample to mean the ancestors of an allele at the neutral locus, as opposed to an individual in the diploid population. Since the transitions of the model can be interpreted as acting on haploids, this should cause no confusion.) The probability of identity will depend on the past history of the population. If we knew that history then we could calculate the identities by iterating backwards in time. We can still make progress if we assume that the population is drawn from a stationary distribution. The Wright–Fisher model described above is just a finite state space Markov chain. The probability of going from $i$ copies of $P$ at time $t$ to $j$ at time $t+1$ is given by

$$P_{ij} = \binom{2N}{j} (p^{**})^j (1-p^{**})^{2N-j}.$$

Provided that it has a nondegenerate stationary distribution $\{\psi_i\}_{i=1}^{2N}$ (which is true provided that the mutation rates $\mu_1$ and $\mu_2$ are strictly positive), then we can reverse the process. The transition probabilities for the backwards in time evolution of the number of type $P$ genomes in the population are given by the prescription

$$\Gamma_{ji} = \frac{\psi_i}{\psi_j} P_{ij}.$$

In following the history of the sample at the neutral locus we must decide whether each individual was associated with a type $P$ or a type $Q$ at the selected locus in the previous generation. This association can change from parent to child as a result of mutation or of recombination. Following through the reproductive process above we see that after the selection and mutation, the proportion of the $P$-population that has arisen by mutation is

$$m_P = \frac{\mu_2 q^*}{(1-\mu_1)p^* + \mu_2 q^*},$$

where we have used the notation $q^* = 1 - p^*$. Similarly, the proportion of the $Q$-gametes that have arisen by mutation is

$$m_Q = \frac{\mu_1 p^*}{(1-\mu_2)q^* + \mu_1 p^*}.$$

Recall that the probability of a recombination event between the selected and neutral sites in one generation is $r$. We need to know the probability, $\tilde{m}_P$, that a neutral locus currently associated with a type $P$ background was associated with a type $Q$ background before the effects of mutation and recombination. First observe that if there is a recombination event, the chance that it is with an individual that is type $Q$ after the mutation step and also with one of type $Q$ before the mutation step is $(1-m_Q)(1-p^{**})$. If there is either no recombination event or a recombination with an individual whose type after the mutation step is $P$, then we require that the type $P$ arose by mutation. Thus, writing $q^{**} = 1 - p^{**}$,

$$\tilde{m}_P = r(1-m_Q)q^{**} + (1-rq^{**})m_P.$$

Similarly,

$$\tilde{m}_Q = rp^{**}(1-m_P) + (1-rp^{**})m_Q.$$

We now have all the information that we require to write down recursions for the quantities of interest at the neutral site. In particular, we are in a position to write down recursive equations for the probability of identity in allelic state for a sample from the neutral locus. (It is this quantity that we shall concentrate on in our numerical examples of Section 7.) We assume that in each generation mutation to a novel allele at the neutral locus occurs with probability $\nu$ and also that the frequencies at the selected locus have reached stationarity. We then write $\{f_{PP,i}, f_{PQ,i}, f_{QQ,i}\}$ for the probabilities of identity in state of a sample of two gametes given that the current number of copies of the $P$ allele in the population is $i$. The subscripts $PP$, $PQ$ and $QQ$ designate the type at the selected locus of the two individuals in our sample. Writing $\{f^*_{PP,i}, f^*_{PQ,i}, f^*_{QQ,i}\}$ for the probability of identity after selection, mutation and recombination at the selected site we have, for

$1 < j < 2N-1$,

$$
\begin{aligned}
f^*_{PP,j} &= \sum_i \Gamma_{ji}(f_{QQ,i}\tilde{m}_P^2 + 2f_{PQ,i}\tilde{m}_P(1-\tilde{m}_P) + f_{PP,i}(1-\tilde{m}_P)^2),\\
f^*_{PQ,j} &= \sum_i \Gamma_{ji}\Big(f_{QQ,i}(1-\tilde{m}_Q)\tilde{m}_P\\
&\qquad\qquad + f_{PQ,i}(\tilde{m}_P\tilde{m}_Q + (1-\tilde{m}_Q)(1-\tilde{m}_P))\\
&\qquad\qquad\qquad\qquad + f_{PP,i}\tilde{m}_Q(1-\tilde{m}_P)\Big),\\
f^*_{QQ,j} &= \sum_i \Gamma_{ji}(f_{QQ,i}(1-\tilde{m}_Q)^2 + 2f_{PQ,i}\tilde{m}_Q(1-\tilde{m}_Q) + f_{PP,i}\tilde{m}_Q^2).
\end{aligned}
\tag{1}
$$

After taking into account random sampling and mutation at the neutral locus, the identities become

$$
\begin{aligned}
f_{PP,j} &= (1-\nu)^2\left(f^*_{PP,j} + \frac{(1-f^*_{PP,j})}{j}\right),\\
f_{PQ,j} &= (1-\nu)^2 f^*_{PQ,j},\\
f_{QQ,j} &= (1-\nu)^2\left(f^*_{QQ,j} + \frac{(1-f^*_{QQ,j})}{2N-j}\right),
\end{aligned}
\tag{2}
$$

with $f_{PP,1} = 1$, $f_{QQ,2N-1} = 1$.

Armed with these equations, numerical calculation of probability of identity now amounts to iterating matrix equations to find a fixed point. However, there is an obstruction to studying this system of equations analytically. Although with strictly positive mutation rates the Wright–Fisher model must have a stationary distribution, it is not known explicitly and consequently, neither are the transition probabilities $\Gamma_{ij}$. In our simulations of Section 7 these are calculated numerically. The numerical estimates show that the Wright–Fisher model is not reversible (that is $\Gamma_{ij}$ does not coincide with $P_{ij}$). In the next section we write down a diffusion approximation for the Moran version of this model. The distinction between the Moran and Wright–Fisher models is that in the Moran model we have overlapping generations. This has the advantage that the frequency of $P$-alleles will then be a generalized birth death process and consequently has a unique invariant measure and this invariant measure is *reversible*. The backward in time transition probabilities for the proportion of type $P$ are then known: they are just given by the forward in time transition probabilities. From Ethier and Kurtz [(1986), Chapter 10, Section 2] we can check that the diffusion approximations for the Wright–Fisher model and that found here for the Moran model are the same. Moreover, since both models are exchangeable, the genealogies of a sample from the population predicted by the two models will coincide in the diffusion limit [see Kingman (1982a)].

To identify the appropriate Moran version of our model, consider again a diploid population, but now evolving according to a continuous time Markov chain so that, in particular, generations overlap. A single step of the chain corresponds to the death of one (diploid) individual and its replacement by another. Such deaths occur at exponential rate $N$ and each individual in the population is equally likely to die. The reproductive step follows the same sequence as before. Writing, as before, $p^*$ and $q^*$ for the proportions of the two types of gamete after the action of selection, and $p^{**}$, $q^{**}$ for the corresponding proportions after both selection and mutation, we have Table 1.

For convenience, we write $(1+s)/(2+s) = (1+S)/2$. Since in the Wright–Fisher model an individual chooses her parents at random, we see that the natural continuous time analogue of our Wright–Fisher model for the evolution of allele frequencies at the selected site is the following version of the Moran model.

Definition 2.1 (The Moran model). The Moran model of a population of size $2N$ is a continuous time Markov chain. At exponential rate $N$, a pair of individuals is chosen at random from the population. One dies and the other reproduces. If the pair chosen consists of one type $P$ and one type $Q$ individual, then the probability that it is the $P$ individual that reproduces is $(1+S)/2$. A type $P$ parent produces a type $P$ offspring with probability $1-\mu_1$, otherwise her offspring is type $Q$. Similarly, a type $Q$ parent has type $Q$ offspring with probability $1-\mu_2$, otherwise her offspring is type $P$.

Remark 2.2. (i) As before, we can take the parameter $s$, and consequently $S$, to be frequency dependent.

(ii) We are assuming that the sampling at each birth/death event is with replacement. This simplifies the expressions for the transition probabilities of the process of allele frequencies at the selected locus, but the price that we pay is that it will somewhat complicate those for the transition probabilities for the sample as we trace backwards in time. Whether we sample with or without replacement will not change the diffusion limit.

Table 1

| **Type to die** | $p^*$ | $p^{**}$ |
|---|---|---|
| $PP$ | 1 | $(1-\mu_1)$ |
| $PQ$ | $(1+s)/(2+s)$ | $(1-\mu_1)p^* + \mu_2 q^*$ |
| $QQ$ | 0 | $\mu_2$ |

In order to keep track of the type at the selected locus of individuals in our sample from the neutral locus, we must also incorporate the effects of recombination. As before, we suppose that at each birth/death event there is a probability $r$ of a recombination event. The type at the selected locus of the offspring of such an event is then inherited not from her parent, but from the individual that died. By this mechanism, as well as by mutation at the selected site, we see migration between the two genetic backgrounds.

Suppose then that we have a sample of individuals from our population and that the type at the selected locus of each individual in our sample is known. We write $n_1^{(1)}(0)$ for the number of individuals in the sample in background $P$ and $n_2^{(1)}(0)$ for the number in background $Q$. We are concerned with the ancestry of the sample. (We superimpose the effects of mutation to a novel type at the neutral locus later.) Thus we write $n_1^{(1)}(t)$ for the number of ancestors associated with type $P$ at the selected locus and $n_2^{(1)}(t)$ for the number of ancestors associated with type $Q$ at time $t$ before the present. We write $p^{(1)}(t)$ for the proportion of the whole population that are type $P$ at that time.

Our final task in this section is to write down the generator of the (backward in time) process $\{(p^{(1)}(t), n_1^{(1)}(t), n_2^{(1)}(t))\}_{t\geq 0}$. We suppose that $p^{(1)}(0)$ is drawn from the stationary distribution of $\{p^{(1)}(t)\}_{t\geq 0}$ and $(n_1^{(1)}(0), n_2^{(1)}(0))$ is arbitrary. The generator will be a very cumbersome object. Mercifully, things will be greatly simplified when we pass to a diffusion approximation.

As we remarked above, the stationary distribution for the Moran model is reversible, so the backwards in time dynamics of the allele frequency, $\{p^{(1)}(t)\}_{t\geq 0}$, are the same as the forwards in time ones described in Definition 2.1. We are going to need the transition probabilities for this process.

LEMMA 2.3. *Consider the probabilities in the jump chain of* $\{p^{(1)}(t)\}_{t\geq 0}$. *Suppose that—looking backward in time—the proportion of the population of type $P$ immediately before and after an arbitrary birth/death event are $p_m$ and $p_{m+1}$, respectively. Writing $P_{p,\tilde{p}} = \mathbb{P}[p_{m+1} = \tilde{p} | p_m = p]$, we have*

$$\begin{aligned} P_{p,p} &= p^2(1-\mu_1) + (1-p)^2(1-\mu_2) \\ &\quad + 2p(1-p)\left(\frac{1+S}{2}\mu_1 + \frac{1-S}{2}\mu_2\right), \end{aligned} \tag{3}$$

$$P_{p,p-1/(2N)} = p\left(p\mu_1 + 2(1-p)\frac{1-S}{2}(1-\mu_2)\right), \tag{4}$$

$$P_{p,p+1/(2N)} = (1-p)\left((1-p)\mu_2 + 2p\frac{1+S}{2}(1-\mu_1)\right). \tag{5}$$

We can now write down the generator of the triple $\{(p^{(1)}(t), n_1^{(1)}(t), n_2^{(1)}(t))\}_{t\geq 0}$. We write

$$p_- = p - \frac{1}{2N}, \qquad p_+ = p + \frac{1}{2N},$$
$$q_- = (1-p_-), \qquad q = (1-p), \qquad q_+ = (1-p_+).$$

LEMMA 2.4. *The generator, $A_{(1)}$, of the process $\{(p^{(1)}(t), n_1^{(1)}(t), n_2^{(1)}(t))\}_{t\geq 0}$ is given by*

$$\begin{aligned}
&A_{(1)} f(p, n_1, n_2) \\
&\quad = N(1-r)c_-\left\{\frac{\binom{n_1}{2}}{\binom{2Np}{2}} p_- q_- (1+S_-)(1-\mu_1) + \frac{n_1 n_2}{4N^2 pq} q_- q \mu_2\right\} \\
&\qquad \times (f(p_-, n_1 - 1, n_2) - f(p, n_1, n_2)) \\
&\qquad + Nc_-\left\{\frac{n_1(2Nq - n_2)}{4N^2 pq} q_- q\mu_2\right. \\
&\qquad\qquad \left. + \frac{rn_1}{2Np} p_- q_- (1+S_-)(1-\mu_1) + \frac{n_1}{2Np} q_- \frac{1}{2N}\mu_2\right\} \\
&\qquad \times (f(p_-, n_1 - 1, n_2 + 1) - f(p, n_1, n_2)) \\
&\qquad + (NP_{p,p-1/(2N)} - R_1 - R_2)(f(p_-, n_1, n_2) - f(p, n_1, n_2)) \\
&\qquad + N(1-r)c_+\left\{\frac{\binom{n_2}{2}}{\binom{2Nq}{2}} p_+ q_+ (1-S_+)(1-\mu_2) + \frac{n_1 n_2}{4N^2 pq} p_+ p \mu_1\right\} \\
&\qquad \times (f(p_+, n_1, n_2 - 1) - f(p, n_1, n_2)) \\
&\qquad + Nc_+\left\{\frac{n_2(2Np - n_1)}{4N^2 pq} p_+ p\mu_1\right. \\
&\qquad\qquad \left. + \frac{rn_2}{2Nq} p_+ q_+ (1-S_+)(1-\mu_2) + \frac{n_2}{2Nq} p_+ \frac{1}{2N}\mu_1\right\} \\
&\qquad \times (f(p_+, n_1 + 1, n_2 - 1) - f(p, n_1, n_2)) \\
&\qquad + (NP_{p,p+1/(2N)} - R_4 - R_5)(f(p_+, n_1, n_2) - f(p, n_1, n_2)) \\
&\qquad + N(1-r)\left\{\frac{\binom{n_1}{2}}{\binom{2Np}{2}} pp_-(1-\mu_1) + \frac{n_1 n_2}{4N^2 pq} pq(1-S)\mu_2\right\} \\
&\qquad \times (f(p, n_1 - 1, n_2) - f(p, n_1, n_2))
\end{aligned}$$

$$
\begin{aligned}
&+ N(1-r)\left\{\frac{\binom{n_2}{2}}{\binom{2Nq}{2}} qq_+(1-\mu_2) + \frac{n_1 n_2}{4N^2 pq} pq(1+S)\mu_1\right\} \\
&\qquad \times (f(p, n_1, n_2 - 1) - f(p, n_1, n_2)) \\
&+ N(1-r)\frac{n_1(2Nq - n_2)}{4N^2 pq} pq(1-S)\mu_2 \\
&\qquad \times (f(p, n_1 - 1, n_2 + 1) - f(p, n_1, n_2)) \\
&+ N(1-r)\frac{n_2(2Np - n_1)}{4N^2 pq} pq(1+S)\mu_1 \\
&\qquad \times (f(p, n_1 + 1, n_2 - 1) - f(p, n_1, n_2)),
\end{aligned}
$$

*where*

$$
S = S(p) = \frac{s(p)}{2 + s(p)}, \qquad S_- = S(p_-), \qquad S_+ = S(p_+),
$$

$$
\begin{aligned}
c_- &= \frac{P_{p,p_-}}{P_{p_-,p}} = \frac{p(p\mu_1 + (1-p)(1-S)(1-\mu_2))}{(1-p_-)((1-p_-)\mu_2 + p_-(1+S_-)(1-\mu_1))}, \\
c_+ &= \frac{P_{p,p_+}}{P_{p_+,p}} = \frac{(1-p)((1-p)\mu_2 + p(1+S)(1-\mu_1))}{p_+(p_+\mu_1 + (1-p_+)(1-S_+)(1-\mu_2))}
\end{aligned}
$$

*and $R_i$ denotes the rate in the $i$th term of the above expression.*

PROOF. Conditional on the changes in $p$, we calculate the probabilities of the possible changes in $(n_1, n_2)$. A backward in time birth/death event corresponds to a forward in time one. To establish the genealogy of our sample, we need to know the role of individuals in the sample in this forward transition. Viewed backward in time the possible transitions of the sample are "migrations," in which the type of the parent of an individual in the sample differs from that of her offspring, and "coalescences," in which two individuals in the sample arise from the splitting of an individual in the previous generation. Although when we pass to the diffusion limit these processes will not happen in a single step (so that two individuals in the sample of different types will not coalesce), here we cannot exclude that possibility.

First observe that exactly three individuals are involved in the (forward in time) birth/death event: the individual that died, the one that split (i.e., gave birth) and her offspring. We write $(i, j, k)$ with $i, j, k \in \{P, Q\}$ for the event that the types of these three individuals are respectively $i, j, k$. Because we have assumed that we are sampling with replacement, the individual that died can coincide with the one that split. In this case, we write $(i, j, k)$ as $(i,k)$. Let us write $\hat{p}_m, \hat{p}_{m+1}$ for the *forwards* in time process immediately before and after an arbitrary birth/death event. From Bayes' rule,

$$\begin{aligned}\mathbb{P}[(i,j,k)|p_{m+1}=\tilde{p}, p_m=p] &= \mathbb{P}[(i,j,k)|\hat{p}_m=\tilde{p}, \hat{p}_{m+1}=p] \\ &= \frac{\mathbb{P}[(i,j,k)\cap\{\hat{p}_{m+1}=p\}|\hat{p}_m=\tilde{p}]}{P_{\tilde{p},p}} \\ &= \frac{\mathbb{P}[(i,j,k)|\hat{p}_m=\tilde{p}]}{P_{\tilde{p},p}}\chi_{(i,j,k),\tilde{p},p},\end{aligned}$$

where $\chi_{(i,j,k),\tilde{p},p}$ is one if the event $(i,j,k)$ results in a (forward in time) change in the proportion of type $P$ from $\tilde{p}$ to $p$ and zero otherwise.

We now use this prescription to calculate all nonzero conditional probabilities of this form. Again we write $p_m, p_{m+1}$ for the proportions in the jump chain backwards in time. First suppose that the gene frequency does not change:

$$\begin{aligned}\mathbb{P}[(P,P,P)|p_m=p, p_{m+1}=p] &= \frac{p(p-1/(2N))(1-\mu_1)}{P_{p,p}}, \\ \mathbb{P}[(P,P)|p_m=p, p_{m+1}=p] &= \frac{(p/(2N))(1-\mu_1)}{P_{p,p}}, \\ \mathbb{P}[(Q,Q,Q)|p_m=p, p_{m+1}=p] &= \frac{(1-p)(1-p-1/(2N))(1-\mu_2)}{P_{p,p}}, \\ \mathbb{P}[(Q,Q)|p_m=p, p_{m+1}=p] &= \frac{(1-p)(1/(2N))(1-\mu_2)}{P_{p,p}}, \\ \mathbb{P}[(P,Q,P)|p_m=p, p_{m+1}=p] &= \frac{p(1-p)(1-S)\mu_2}{P_{p,p}}, \\ \mathbb{P}[(Q,P,Q)|p_m=p, p_{m+1}=p] &= \frac{p(1-p)(1+S)\mu_1}{P_{p,p}}.\end{aligned}$$

When a type $P$ individual is lost (looking backward in time) we have

$$\mathbb{P}\left[(Q,Q,P)\Big|p_m=p, p_{m+1}=p-\frac{1}{2N}\right] = \frac{(1-p+1/(2N))(1-p)\mu_2}{P_{p-1/(2N),p}},$$

$$\mathbb{P}\left[(Q,P)\Big|p_m=p, p_{m+1}=p-\frac{1}{2N}\right] = \frac{(1-p+1/(2N))1/(2N)\mu_2}{P_{p-1/(2N),p}},$$

$$\begin{aligned}&\mathbb{P}\left[(Q,P,P)\Big|p_m=p, p_{m+1}=p-\frac{1}{2N}\right] \\ &\quad= \frac{(p-1/(2N))(1-p+1/(2N))(1+S_-)(1-\mu_1)}{P_{p-1/(2N),p}}.\end{aligned}$$

Finally, if a type one individual is gained (looking backward in time) we have

$$\mathbb{P}\left[(P,P,Q)\Big|p_m = p, p_{m+1} = p + \frac{1}{2N}\right] = \frac{(p+1/(2N))p\mu_1}{P_{p+1/(2N),p}},$$

$$\mathbb{P}\left[(P,Q)\Big|p_m = p, p_{m+1} = p + \frac{1}{2N}\right] = \frac{(p+1/(2N))(1/(2N))\mu_1}{P_{p+1/(2N),p}},$$

$$\mathbb{P}\left[(P,Q,Q)\Big|p_m = p, p_{m+1} = p + \frac{1}{2N}\right]$$
$$= \frac{(p+1/(2N))(1-p-1/(2N))(1-S_+)(1-\mu_2)}{P_{p+1/(2N),p}}.$$

All that remains is to establish the probability that these birth/death events involved individuals in the sample.

First we consider coalescence. A coalescence of two individuals associated with type $P$ occurs in the sample if the parent and offspring are both associated with type $P$ and both form part of the sample and there was no recombination. Thus if $p_m = p$, conditional on $(P,P,P)$ or $(Q,P,P)$, this happens with probability $(1-r)\binom{n_1}{2}/\binom{2Np}{2}$. Similarly for a coalescence of two individuals in the sample associated with type $Q$. For individuals associated with type $P$ and $Q$ from the sample to coalesce requires an event of the form $(i,Q,P)$ or $(i,P,Q)$ (corresponding to parent and offspring having different type) and conditional on one of these events happening, has probability $(1-r)\frac{n_1 n_2}{2Np(2N-2Np)}$.

Now we consider "migration." An individual in the sample can "migrate" from one background to the other as a result of mutation or recombination. In either case she must be the offspring of a birth/death event. If the parent and the individual that die have different types, then a recombination combined with a mutation does not lead to a change in background. Thus conditional on $p_m = p$ and a birth/death event that involved mutation from type $Q$ to $P$ forward in time, an individual in our sample will migrate from type $P$ to $Q$ (backward in time) with probability $(1-r)\frac{n_1(2N-2Np-n_2)}{2Np(2N-2Np)}$ if the individual that dies and the individual that split are different types, $\frac{n_1(2N-2Np-n_2)}{2Np(2N-2Np)}$ if they are the same type but different individuals and $n_1/(2Np)$ if the individual that dies is the parent. Similarly, conditional on $p_m = p$ and a birth/death event involving a mutation from $P$ to $Q$, an individual in our sample migrates from type $Q$ to $P$ (backward in time) with probability $(1-r)\frac{n_2(2Np-n_1)}{2Np(2N-2Np)}$ if the event is $(Q,P,Q)$, with probability $\frac{n_2(2Np-n_1)}{2Np(2N-2Np)}$ if the event is $(P,P,Q)$ and with probability $\frac{n_2}{2N-2Np}$ if the event is $(P,Q)$.

If there is no mutation, an individual in the sample can still change type due to recombination: conditional upon $(Q,P,P)$ with probability $\frac{rn_1}{2Np}$, a member of the sample migrates from background $P$ to $Q$ and conditional on $(P,Q,Q)$ with probability $\frac{rn_2}{2N-2Np}$, an individual of the sample migrates from background $Q$ to background $P$.

Finally, recalling that events in the jump chain take place at an exponential rate $N$, we obtain the claimed expression. $\square$

**3. The generator of the diffusion approximation.** In this section we identify the generator of the diffusion approximation corresponding to the backward in time model of Section 2. Existence of a corresponding stochastic process and convergence to the limit is deferred to the following sections, but our proofs will require the following assumption on the selection coefficient.

ASSUMPTION. The selection coefficient, $s:[0,1]\to\mathbb{R}$ is a Lipschitz continuous function.

As usual, we speed up time by a factor of diploid population size, $N$, and correspondingly scale down the parameters in the model by the same factor. Thus $\mu_i\mapsto\mu_i/N$ and $r\mapsto r/N$. The selection coefficient, $s$, is also scaled by $N$. Notice that

$$\frac{1+s/N}{2+s/N}=\frac{1}{2}\left(1+\frac{s}{2N}\right)+o\left(\frac{1}{N}\right),$$

so that at the $N$th stage of the rescaling $S=s/(2N)+o(1/N)$.

We write $A_{(N)}$ for the generator of Lemma 2.4 with parameters scaled in this way.

LEMMA 3.1. *Let* $E=[0,1]\times\{1,\dots,n_1(0)+n_2(0)\}\times\{1,\dots,n_1(0)+n_2(0)\}$ *and suppose that* $f(p,n_1,n_2):E\to\mathbb{R}$ *is twice continuously differentiable with respect to* $p$. *Then for* $0<p<1$,

$$A_{(N)}f(p,n_1,n_2)\to Af(p,n_1,n_2)\qquad \text{as } N\to\infty,$$

*where*

$$\begin{aligned}
&Af(p,n_1,n_2)\\
&\quad=\frac{1}{2p}\binom{n_1}{2}(f(p,n_1-1,n_2)-f(p,n_1,n_2)) &(6)\\
&\qquad+\frac{1}{2q}\binom{n_2}{2}(f(p,n_1,n_2-1)-f(p,n_1,n_2)) &(7)\\
&\qquad+\frac{p}{q}\mu_1\frac{n_2}{2}(f(p,n_1+1,n_2-1)-f(p,n_1,n_2)) &(8)
\end{aligned}$$

$$+\frac{q}{p}\mu_2\frac{n_1}{2}(f(p,n_1-1,n_2+1)-f(p,n_1,n_2)) \tag{9}$$

$$+r\frac{n_2p}{2}(f(p,n_1+1,n_2-1)-f(p,n_1,n_2)) \tag{10}$$

$$+r\frac{n_1q}{2}(f(p,n_1-1,n_2+1)-f(p,n_1,n_2)) \tag{11}$$

$$+(-\mu_1p+\mu_2q+spq)\frac{1}{2}f'(p,n_1,n_2)+\frac{1}{4}pqf''(p,n_1,n_2), \tag{12}$$

*and "′" denotes differentiation with respect to p.*

REMARK 3.2. Note that $\mu_1$, the mutation rate for $P\to Q$ forward in time, is involved in jumps $Q\to P$ in this backward in time generator (analogous for $\mu_2$).

PROOF OF LEMMA 3.1. Fix $p>0$ and consider the generator $A_{(N)}$, which is the generator $A_{(1)}$ of Lemma 2.4 with the parameters scaled as above. Notice that $c_+$ and $c_-$ tend to one as $N\to\infty$ and that the effect of speeding up time is to multiply the whole generator by a further factor of $N$.

Like $A_{(1)}$ the generator $A_{(N)}$ consists of ten terms corresponding to all possible events. The first and seventh term sum to give

$$\frac{1}{2p}\binom{n_1}{2}(f(p,n_1-1,n_2)-f(p,n_1,n_2))+O\left(\frac{1}{N}\right).$$

The fourth and eighth terms sum to give

$$\frac{1}{2q}\binom{n_1}{2}(f(p,n_1,n_2-1)-f(p,n_1,n_2))+O\left(\frac{1}{N}\right).$$

Ignoring the part arising from recombination, the fifth and tenth terms sum to give

$$\frac{p}{q}\mu_1\frac{n_2}{2}(f(p,n_1+1,n_2-1)-f(p,n_1,n_2))+O\left(\frac{1}{N}\right).$$

Ignoring the part arising from recombination, the second and ninth terms sum to give

$$\frac{q}{p}\mu_2\frac{n_1}{2}(f(p,n_1-1,n_2+1)-f(p,n_1,n_2))+O\left(\frac{1}{N}\right).$$

The contribution to the second and fifth terms from recombination sum to give

$$r\frac{n_2p}{2}(f(p,n_1+1,n_2-1)-f(p,n_1,n_2)) \\ +r\frac{n_1q}{2}(f(p,n_1-1,n_2+1)-f(p,n_1,n_2))+O\left(\frac{1}{N}\right).$$

This leaves the third and sixth terms. First observe that $NR_i$ for $i = 1, 2, 4, 5$, are all $O(1)$, whereas $N^2 P_{p,p-1/(2N)}$ and $N^2 P_{p,p+1/(2N)}$ will both be $O(N^2)$. Using smoothness of $f$ as a function of $p$, we expand $f(p_-, n_1, n_2)$ and $f(p_+, n_1, n_2)$ in a Taylor series about $(p, n_1, n_2)$. The third and sixth terms then become

$$(N^2 P_{p,p-1/(2N)} - NR_1 - NR_2)\left(-\frac{1}{2N} f'(p, n_1, n_2) + \frac{1}{8N^2} f''(p, n_1, n_2)\right)$$
$$+ (N^2 P_{p,p+1/(2N)} - NR_4 - NR_5)\left(\frac{1}{2N} f'(p, n_1, n_2) + \frac{1}{8N^2} f''(p, n_1, n_2)\right)$$
$$+ O\left(\frac{1}{N}\right),$$

where we have again used "$'$" to denote differentiation with respect to $p$. Substituting from Lemma 2.3 we obtain

$$N^2 p\left(p\frac{\mu_1}{N} + q\left(1 - \frac{s}{2N}\right)\left(1 - \frac{\mu_2}{N}\right)\right)$$
$$\times \left(-\frac{1}{2N} f'(p, n_1, n_2) + \frac{1}{8N^2} f''(p, n_1, n_2)\right)$$
$$+ N^2 q\left(q\frac{\mu_2}{N} + p\left(1 + \frac{s}{2N}\right)\left(1 - \frac{\mu_1}{N}\right)\right)$$
$$\times \left(\frac{1}{2N} f'(p, n_1, n_2) + \frac{1}{8N^2} f''(p, n_1, n_2)\right) + O\left(\frac{1}{N}\right),$$

which reduces to

$$\left(-p^2\mu_1 + \frac{pqs}{2} + \mu_2 pq\right)\frac{1}{2} f'(p, n_1, n_2)$$
$$+ \left(q^2\mu_2 + \frac{pqs}{2} - \mu_1 pq\right)\frac{1}{2} f'(p, n_1, n_2) + \frac{1}{4} pq f''(p, n_1, n_2) + O\left(\frac{1}{N}\right)$$
$$= (-\mu_1 p + \mu_2 q + spq)\frac{1}{2} f'(p, n_1, n_2) + \frac{1}{4} pq f''(p, n_1, n_2) + O\left(\frac{1}{N}\right).$$

Letting $N \to \infty$ in the above expressions completes the proof. □

REMARK 3.3. There are no surprises in the form of this generator. If we think of our population as subdivided into two demes according to the type at the selected site, then we should expect the genealogy of the sample to be given by a structured coalescent. Thus the terms (6) and (7) correspond to the coalescence of individuals in the same deme, which happens at a rate inversely proportional to the population size within that deme. The terms (8) and (9) reflect migration between demes as a result of mutation. The rates must be scaled by the ratio of the population sizes in the different demes, just as in the structured coalescent. The terms (10) and (11) reflect migration due to recombination. Evidently these rates must be proportional to the proportion of the population that is of the opposite type. (Recombining with an individual of one's own type has no net effect.) Finally the term (12) is simply the generator of the diffusion approximation to our Moran model for allele frequencies.

**4. Existence of the diffusion approximation.** It is not immediately obvious that there should be a stochastic process with generator given by (6)–(12). The immediate problem is that the coalescence and migration rates for the sample become unbounded as the allele frequency $p$ tends to zero or one. This means that, in principle, we could see an infinite number of jumps in finite time. However, what we shall see is that this does not happen because the process jumps away from the "bad region" in a finite number of jumps.

First let us define "bad" (or rather "good" regions) for the process. Evidently we want to keep away from regions where $p$ is small and $n_1 \neq 0$ or where $q$ is small and $n_2 \neq 0$. We therefore define

$$\begin{aligned} U^{(k)} &= \left[0, \frac{1}{k}\right] \times \{0\} \times \{0, 1, \ldots, n_1(0) + n_2(0)\} \\ &\cup \left[1 - \frac{1}{k}, 1\right] \times \{0, 1, \ldots, n_1(0) + n_2(0)\} \times \{0\} \\ &\cup \left(\frac{1}{k}, 1 - \frac{1}{k}\right) \times \{0, 1, 2, \ldots, n_1(0) + n_2(0)\}\{0, 1, 2, \ldots, n_1(0) + n_2(0)\}. \end{aligned}$$

Notice that the sets $U^{(k)}$ are open subsets of the state space $E$.

The first, straightforward, task is to show that the process exists until its exit time from $U^{(k)}$ for each $k$. In an obvious notation, we write $A_n$ for the portion of the generator corresponding to the terms (6)–(11) and $A_p$ for the portion corresponding to the generator of the allele frequencies, namely (12).

LEMMA 4.1. *We write $C^2(E)$ for bounded functions $f: E \to \mathbb{R}$ which are twice continuously differentiable with respect to $p$. Define*

$$A_p^{(k)} f(p, n_1, n_2) = \chi_{U^{(k)}}((p, n_1, n_2)) A_p f(p, n_1, n_2),$$

$$A_n^{(k)} f(p, n_1, n_2) = \chi_{U^{(k)}}((p, n_1, n_2)) A_n f(p, n_1, n_2).$$

*Then the closure of*

$$\{(f, A^{(k)} f) : f \in C^2(E)\} \equiv \{(f, A_p^{(k)} f + A_n^{(k)} f) : f \in C^2(E)\}$$

*generates a Feller semigroup.*

PROOF. This is standard. First consider $A_p^{(k)}$ applied to functions $f \in C^2([0,1])$. Then from, for example, Ethier and Kurtz [(1986), Chapter 8, Theorem 2.8], the closure of $\{(f, A_p f) : f \in C^2([0,1])\}$ generates a Feller semigroup on the continuous functions on $[0,1]$. (It is here that we have used that the selection coefficient, $s$, is a *Lipschitz* continuous function.) We also know that, for fixed $p \in [0,1]$, $A_n^{(k)}$ generates a Feller semigroup on continuous functions on $\{0,1,\dots,n_1(0)+n_2(0)\} \times \{0,1,\dots,n_1(0)+n_2(0)\}$ (since the state space is finite and the jump rates are all bounded). Evidently both generators can be regarded as acting on $E$ and they are Feller generators on continuous functions on $E$.

Now observe that

$$\|A_n^{(k)} f\|_\infty \le (n_1(0)^2 + n_2(0)^2 + (2k+r)(n_1(0)+n_2(0)))\|f\|_\infty$$

and so $A^{(k)}$ is a bounded perturbation of $A_p^{(k)}$ and hence also generates a strongly continuous contraction semigroup [see Ethier and Kurtz (1986), Chapter 1, Section 7]. That the resulting semigroup is positive and conservative is an easy consequence of the Trotter product formula and the proof is complete. □

We can think of the processes constructed in Lemma 4.1 as solutions to a *stopped* martingale problem. Let $\tau_k = \inf\{t \ge 0 : (p(t), n_1(t), n_2(t)) \notin U^{(k)}\}$. Then, for $f$ in the domain of $A$,

$$\begin{aligned} &f(p(t), n_1(t), n_2(t)) - \int_0^{t\wedge\tau_k} A^{(k)} f(p(s), n_1(s), n_2(s))\, ds \\ &\quad = f(p(t), n_1(t), n_2(t)) - \int_0^{t\wedge\tau_k} A f(p(s), n_1(s), n_2(s))\, ds \end{aligned} \tag{13}$$

is a martingale. We remark that this stopped martingale problem is well posed since it is associated with a Feller generator [see Ethier and Kurtz (1986), Chapter 4, Theorem 4.1].

To establish existence of the process corresponding to the generator $A$ on the whole of $E$, we shall use the following result [Ethier and Kurtz (1986), Chapter 4, Theorem 6.3].

THEOREM 4.2. *Let $(E, d)$ be a complete and separable metric space and let $A \subset \overline{C}(E) \times B(E)$. Let $U_1 \subset U_2 \subset \cdots$ be open subsets of $E$. Fix $\nu \in \mathcal{P}(E)$ and suppose that for each $k$ there exists a unique solution $X_k$ of the stopped martingale problem for $(A, \nu, U_k)$ with sample paths in $D_E[0,\infty)$. Setting*

$$\tau_k = \inf\{t : X_k(t) \notin U_k \text{ or } X_k(t-) \notin U_k\},$$

*suppose that for each $t > 0$,*

$$\lim_{k\to\infty} \mathbb{P}\{\tau_k \le t\} = 0. \tag{14}$$

*Then there exists a unique solution of the $D_E[0,\infty)$ martingale problem for $(A,\nu)$.*

Here $\overline{C}(E)$ denotes bounded continuous functions on $E$, $B(E)$ denotes bounded Borel measurable functions on $E$, $D_E[0,\infty)$ is càdlàg paths in $E$ and $\mathcal{P}(E)$ is probability measures on $E$. The generator $A$ is to be thought of as $\{(f,Af)\}$ for $f$ in a suitable class. [In our case $A$ should be thought of as the closure of $\{(f,Af) : f \in C^2(E)\}$.] The probability measure $\nu$ specifies the initial distribution of our process.

Our task then is to show that if we take $U_k = U^{(k)}$, then the condition (14) is satisfied.

PROPOSITION 4.3. With $\tau_k$ as above, for any fixed $t > 0$,

$$\lim_{k\to\infty} \mathbb{P}\{\tau_k \le t\} = 0.$$

Of course, if the boundaries are inaccessible for the process $\{p(t)\}_{t\ge 0}$, which, as we show as part of Lemma 4.4, is true provided $\mu_i \ge 1/2$ for $i = 1,2$, then there is no problem. The difficulty is to check that this is still true under the much weaker (and biologically more realistic) condition that we are assuming here, that $\mu_i > 0$ for $i = 1,2$. The key is the following lemma.

LEMMA 4.4. (i) *Suppose that $\mu_2 \ge 1/2$ (resp. $\mu_1 \ge 1/2$). Then 0 (resp. 1) is an inaccessible boundary for the process $\{p(t)\}_{t\ge 0}$. If $\mu_2 \in (0,1/2)$ [resp. $\mu_1 \in (0,1/2)$], then it is accessible.*

(ii) *Suppose that $\mu_2 < 1/2$ (resp. $\mu_1 < 1/2$). Then for any fixed value of $p(0) \in (0,1)$ and any $K > 0$, writing $\tau_x(a)$ for the first hitting time of $a$ by the process $\{p(t)\}_{t\ge 0}$ given that $p(0) = x$, we have*

$$\lim_{k\to\infty} \mathbb{P}\left[\int_0^{\tau_{p(0)}(1/k)} \frac{1}{p(s)}\,ds > K\right] = 1, \tag{15}$$

*respectively,*

$$\lim_{k\to\infty} \mathbb{P}\left[\int_0^{\tau_{p(0)}(1-1/k)} \frac{1}{1-p(s)}\,ds > K\right] = 1.$$

REMARK 4.5. Equation (15) is not the strongest statement that we could make about the divergence of the integral, but it is the form that we require in the proof of Proposition 4.3.

PROOF OF LEMMA 4.4. Recall that for a one-dimensional diffusion process on the interval $[0,1]$ with generator

$$L = \frac{1}{2}a(x)\frac{d^2}{dx^2} + b(x)\frac{d}{dx},$$

the scale, $n(x)$, and speed, $m(x)$, are defined for $x \in [0,1]$ by

$$n(x) = \int_c^x \exp\biggl(-\int_c^y \frac{2b(z)}{a(z)}\,dz\biggr)\,dy,$$
$$m(x) = \int_c^x \frac{2}{a(y)} \exp\biggl(\int_c^y \frac{2b(z)}{a(z)}\,dz\biggr)\,dy,$$

where $c \in (0,1)$ is fixed arbitrarily. According to Feller's boundary classification, a boundary point $e$ is accessible or inaccessible according as

$$u(e) \triangleq \int_c^e m(x)\,dn(x)$$

is finite or infinite.

For the process of allele frequencies, $\{p(t)\}_{t\geq 0}$, we have

$$a(x) = \tfrac{1}{2}x(1-x), \qquad b(x) = \tfrac{1}{2}(s(x)x(1-x) - \mu_1 x + \mu_2(1-x)).$$

Substituting gives

$$n(x) = \int_c^x \exp\biggl(-\int_c^y 2s(z)\,dz\biggr) y^{-2\mu_2}(1-y)^{-2\mu_1}\,dy,$$
$$m(x) = \int_c^x \frac{4}{y(1-y)} \exp\biggl(\int_c^y 2s(z)\,dz\biggr) y^{2\mu_2}(1-y)^{2\mu_1}\,dy.$$

Since $s(x)$ is bounded and continuous and $\mu_i > 0$ for $i = 1,2$, we have $n(x) \sim \int_{x+} y^{-2\mu_2}\,dy$ as $x \downarrow 0$ which is bounded or unbounded according as $\mu_2 < 1/2$ or $\mu_2 \geq 1/2$. The symmetrical argument applied to the boundary point 1 completes the proof of Part 1 of the lemma.

Now suppose that $\mu_2 < 1/2$. Since $p(0)$ is arbitrary, for the remainder of this proof we shall suppress it in our notation and write $\tau(x)$ for the first hitting time of $x$.

First we convert the process $p(t)$ to natural scale. That is, we study the process $Y(t) \triangleq n(p(t))$. Now

$$Y(t) = Y(0) + \int_0^t \sigma(Y(s))\,dW_s,$$

where $\{W_t\}_{t\geq 0}$ is a standard Brownian motion started from $n(p(0))$ and (using "$'$" to denote differentiation)

$$\sigma(y) = n'(n^{-1}(y))\sqrt{\tfrac{1}{2}n^{-1}(y)(1-n^{-1}(y))}.$$

Without loss of generality, since scale is defined only up to translation, we may assume that $n(0) = 0$ and then as $x \downarrow 0$ our calculations above show that $n(x) \sim x^{1-2\mu_2}$. The quantity that we are interested in is

$$\int_0^{\tau(1/k)} \frac{1}{n^{-1}(Y(s))}\,ds.$$

Writing $Y(t) = W(\gamma(t))$ where $\gamma$ is defined by

$$\int_0^{\gamma(t)} \frac{1}{\sigma(W_s)}\, ds = t,$$

and substituting $r = \gamma(s)$ in the integral, we obtain

$$\int_0^{\tau(1/k)} \frac{1}{p(s)}\, ds = \int_0^{\gamma(\tau(1/k))} \frac{1}{n^{-1}(W_r)} \frac{1}{\sigma(W_r)}\, dr.$$

Now observe that

$$n^{-1}(x) \sim x^{1/(1-2\mu_2)}, \qquad \sigma(x) \sim x^{-2\mu_2/(1-2\mu_2)} x^{1/(2(1-2\mu_2))} \qquad \text{as } x \downarrow 0,$$

and so the behavior of (15) is determined by

$$\int_0^{\gamma(\tau(1/k))} \frac{1}{W_r^\alpha}\, dr, \tag{16}$$

where

$$\alpha = \frac{1}{1-2\mu_2} - \frac{2\mu_2}{1-2\mu_2} + \frac{1}{2(1-2\mu_2)} = 1 + \frac{1}{2(1-2\mu_2)}.$$

Notice that $W(\gamma(\tau(1/k)) = n(1/k)$ implies that $\gamma(\tau(1/k)) \sim \tau^W(k^{-(1-2\mu_2)})$, where $\tau^W(x)$ is the first hitting time of $x$ by the Brownian motion $\{W_t\}_{t\geq 0}$. Finally, since $\mu_2 < \frac{1}{2}$, $\alpha > 1$ and $\lim_{k\to\infty} \gamma(\tau(1/k)) = \tau^W(0)$, we deduce (15). □

Proof of Proposition 4.3. For $a \in [0,1]$, we retain the notation $\tau_x(a)$ for the first hitting time of $a$ by the process $\{p(t)\}_{t\geq 0}$ given that $p(0) = x$.

Suppose that 0 is an inaccessible boundary for the process $\{p(t)\}_{t\geq 0}$. Then we have $\mathbb{P}[\tau_{1/k}(0) \leq t] \to 0$ as $k \to \infty$. A fortiori, the probability that $\{(p(t), n_1(t), n_2(t))\}_{t\geq 0}$ hits $\{\frac{1}{k}\} \times \{i\} \times \{n_1\}$ for $i \neq 0$ before time $t$ tends to zero as $k$ tends to infinity. Similarly, if 1 is inaccessible for the process of allele frequencies, the probability of hitting $\{1 - \frac{1}{k}\} \times \{n_1\} \times \{i\}$ for nonzero $i$ tends to zero as $k \to \infty$.

We concentrate on the case when 0 is an *accessible* boundary for the process $\{p(t)\}_{t\geq 0}$, that is $\mu_2 < 1/2$.

The idea of the proof is simple. If $p$ approaches zero, then (15) implies that the probability of jumping from type $P$ to type $Q$ for individuals in our sample tends to one. On the other hand the jump rate from type $Q$ to type $P$ tends to zero. Therefore, with very high probability (tending to one), if the allele frequency is less than $1/k$, we do not see a type $P$ individual in (the ancestors of) our sample. Combined, if 1 is also accessible, with the symmetric argument as the allele frequency increases to one, we see that the probability of hitting the boundary of $U^{(k)}$ in finite time converges to zero as $k \to \infty$.

We now make this argument more precise. Suppose that $t \geq 0$ is fixed. We use the abbreviation $n(0) = n_1(0) + n_2(0)$. This is an upper bound on the number of individuals in the sample at all times. Now fix $\delta > 0$. For each $p \in (0,1)$, define

$$\lambda_1(p) = n(0)\left(\frac{p\mu_1}{2(1-p)} + \frac{rp}{2}\right)$$

and

$$\lambda_2(p) = \frac{1-p}{p}\frac{\mu_2}{2}.$$

Notice that $\lambda_2$ provides a lower bound for the rate at which individuals in the sample jump away from state $P$ (provided that $n_1 \neq 0$), whereas $\lambda_1$ provides an upper bound for the rate at which they arrive.

Let $T_\varepsilon$ be an exponentially distributed random variable with rate $\lambda_1(\varepsilon)$.

Recall that for $0 \leq a < x < b \leq 1$,

$$\mathbb{P}[\tau_x(a) < \tau_x(b)] = \frac{n(b) - n(x)}{n(b) - n(a)},$$

where the scale function, $n(x)$ was defined in Lemma 4.4. Substituting, we see that by choosing $N$ large enough, we can arrange that

$$\mathbb{P}[\tau_\varepsilon(0) < \tau_\varepsilon(N\varepsilon)] > 1 - \frac{\delta}{8} \qquad \text{uniformly in } \varepsilon < \frac{1}{N}. \tag{17}$$

By choosing $\varepsilon$ to be still smaller if necessary, we arrange that

$$\mathbb{P}[T_{N\varepsilon} > t] > 1 - \frac{\delta}{8}. \tag{18}$$

Now let $X$ be a Poisson random variable with mean $K$. Choose $K$ large enough that

$$\mathbb{P}[X > n(0)] > 1 - \frac{\delta}{8}. \tag{19}$$

Finally suppose that $p(0) \geq \varepsilon$ and (with $\varepsilon$ fixed) using Lemma 4.4 choose $k_0$ large enough that for $k > k_0$,

$$\mathbb{P}\left[\int_{\tau_{p(0)}(\varepsilon)}^{\tau_{p(0)}(1/k)} \lambda_2(p(s))\, ds > K\right] > 1 - \frac{\delta}{8}. \tag{20}$$

Now consider the first time $\tau_{p(0)}(1/k)$ that the process of allele frequencies hits $1/k$. We want to estimate the probability that $n_1(\tau_{p(0)}(1/k)) = 0$. We combine the above estimates as follows. We use equation (17) to restrict our attention to the event that between first hitting $\varepsilon$ and first hitting $1/k$ we always have $p < N\varepsilon$. Equation (18) then allows us to ignore the possibility that there are any jumps of individuals in the sample *into* state $P$ in this time interval. Equation (20) ensures that the rate of jumps out of state $P$ is at least $K$ (provided there are any individuals to jump) and finally equation (19) ensures that all individuals do indeed jump out of state $P$ before the process hits $1/k$. Thus, with probability at least $1-\delta/2$ when the $p$ process hits $1/k$, $n_1 = 0$.

Started from $p = 1/k$ and $n_1 = 0$, we now let the process run until the first time $\mathcal{T}$ that $n_1 \neq 0$. Evidently this is smallest if $p = 1/k$, $n_2 = n(0)$. Moreover, since the rate at which individuals jump into state $P$ increases as $p$ increases, provided that $\varepsilon < 1/2$, this time is stochastically greater than $\mathcal{T}'$, the first time started from $(1/2, 0, n(0))$ that $n_1 \neq 0$. The distribution of $\mathcal{T}'$ is independent of $\varepsilon$ and $\delta$. We want to apply the above argument once again to see that the next time the allele frequency hits $1/k$, the probability that $n_1 = 0$ is at least $1-\delta/2$. The only twist is that we may need to choose $\varepsilon$ smaller still to ensure that the probability that at the time $\mathcal{T}$ the process $p$ is greater than $\varepsilon$ is at least $1-\delta/2$. Notice that this last estimate can be obtained uniformly in $k > 1/\varepsilon$ by choosing $\varepsilon$ so that if we start from $(0, 0, n(0))$, then $p(\mathcal{T}) > \varepsilon$ with probability $1-\delta/2$.

Now we are essentially done. The process hits the boundary $p = 1/k$ with $n_1 \neq 0$ only after a geometric number of hits of $p$ on $1/k$. The success probability for this geometric random variable is at most $\delta$. Each "failure" accrues an additional waiting time bounded below by an independent copy of $\mathcal{T}'$. Since $\delta$ was arbitrary, the proof is complete. □

**5. Convergence.** Having established existence of our candidate diffusion approximation, we now turn to proving that the rescaled processes of Section 3 actually *converge* to this limit. That is, we prove the following theorem.

THEOREM 5.1. *The processes* $\{(p^{(N)}(t), n_1^{(N)}(t), n_2^{(N)}(t))\}_{t\geq 0}$ *corresponding to the generators* $A_{(N)}$ *of Lemma 3.1 converge weakly in* $D_E[0,\infty)$ *as* $N \to \infty$ *to the process* $\{(p(t), n_1(t), n_2(t))\}_{t\geq 0}$ *generated by* $A$.

The main tool in the proof will be the following result which is a special case of Ethier and Kurtz [(1986), Chapter 4, Corollary 8.7].

THEOREM 5.2. *Suppose that* $(E, d)$ *is a complete separable metric space. Let* $A$ *be a Feller generator on* $E$ *corresponding to the Markov process* $X$. *For each* $N \geq 1$, *let* $X^{(N)}$ *be progressively measurable* $E$*-valued processes with full generators* $\hat{A}^{(N)}$ *and such that* $X^{(N)}(0)$ *converges weakly to* $X(0)$ *as* $N \to \infty$. *Suppose that* $\overline{\mathcal{D}(A)}$ *separates points. Suppose further that the compact containment condition holds for* $\{X^{(N)}\}_{N\geq 1}$. *That is, for every* $\varepsilon > 0$ *and every* $T > 0$ *there exists a compact set* $\Gamma_{\varepsilon,T} \subseteq E$ *for which*

$$\inf_N \mathbb{P}[X^{(N)}(t) \in \Gamma_{\varepsilon,T} \text{ for } 0 \leq t \leq T] \geq 1 - \varepsilon.$$

*Suppose that for each* $(f,g) \in A$ *and* $T > 0$ *there exist* $(f^{(N)}, g^{(N)}) \in \hat{A}^{(N)}$ *and* $G^{(N)} \subset E$ *such that*

$$\lim_{N\to\infty} \mathbb{P}[X^{(N)}(t) \in G^{(N)}, 0 \le t \le T] = 1, \tag{21}$$

$\sup_N \|f^{(N)}\|_\infty < \infty$ *and*

$$\lim_{N\to\infty} \sup_{x\in G^{(N)}} |f(x) - f^{(N)}(x)| = 0 = \lim_{N\to\infty} \sup_{x\in G^{(N)}} |g(x) - g^{(N)}(x)|. \tag{22}$$

*Then* $X^{(N)}$ *converges weakly to* $X$ *as* $N \to \infty$.

In fact, this result would be sufficient to allow us to prove the result in one fell swoop, but in view of the work of Section 4, it is convenient to proceed in two stages. First, taking the sets $G^{(N)} = E$, we prove convergence of the *stopped* processes, $X_{(N,k)}$, generated by $A_{(N,k)} \equiv \chi_{U^{(k)}} A_{(N)}$, to the process $X_k$ for each $k$. The following lemma, which is a straightforward adaptation of Lemma 11.1.1 of Stroock and Varadhan (1979) then completes the proof.

LEMMA 5.3. *Let* $\{\mathbb{P}^{(N)}\}_{N\ge1}$ *be a sequence of probability measures on the space* $D_E[0,\infty)$ *and suppose that* $T^{(k)}$ *is a nondecreasing sequence of stopping times (with respect to the natural filtration) increasing to infinity almost surely. For each* $k \ge 1$, *let* $\{\mathbb{P}^{(N,k)}\}_{N\ge1}$ *be a relatively compact sequence of probability measures such that* $\mathbb{P}^{(N,k)}$ *is equal to* $\mathbb{P}^{(N)}$ *on* $\mathcal{F}_{T^{(k)}}$.

*If the probability measure* $\mathbb{P}$ *has the property that, for any* $k \ge 1$, *any limit point of* $\{\mathbb{P}^{(N,k)}\}_{N\ge1}$ *agrees with* $\mathbb{P}$ *on* $\mathcal{F}_{T^{(k)}}$, *then* $\mathbb{P}^{(N)}$ *converges to* $\mathbb{P}$ *as* $N \to \infty$.

PROOF OF THEOREM 5.1. First we fix $k \ge 1$ and consider the sequence of stopped processes $\{X_{(N,k)}\}_{N\ge1}$. The compact containment condition of Theorem 5.2 is automatically satisfied since $E$ is compact. We take the sets $G_n = E$ in condition (21). From Section 3 we see that we can take $f^{(N)} = f$ in condition (22) and convergence of the stopped processes is proved.

Combining Proposition 4.3 with Lemma 5.3, the proof is complete. □

**6. Differential equations for the identities.** We now return to the problem of calculating the probability of identity in allelic state at the neutral locus for a sample whose types at the selected locus are known. For simplicity, we consider the case of a sample of size two, but see Remark 6.2.

Our approach is to use the diffusion approximation to write down a coupled system of ordinary differential equations for the probability of identity, indexed by the state at the selected locus of the sample, that is $PP$, $PQ$ and $QQ$. Thus $f_{PP}$ will denote the probability of identity given that the two individuals in our sample are both of type $P$. These quantities will be compared to the predictions of equations (1) and (2) in Section 7.

If the current allele frequency at the selected locus is $p$ (assumed as always to have reached stationarity), then write $F_{PP}(t,p)$ for the probability that at time $t$ the (backwards in time) process $\{(p(t), n_1(t), n_2(t))\}_{t\geq 0}$ is in $[0,1]\times\{0\}\times\{1\}\cup[0,1]\times\{1\}\times\{0\}$ given that at time zero $n_1(0)=2$, $n_2(0)=0$. Similarly, define $F_{PQ}(t,p)$ and $F_{QQ}(t,p)$. We assume, as always, that $\{p(0)\}_{t\geq 0}$ is drawn from the (reversible) stationary distribution for the process $\{p(t)\}_{t\geq 0}$.

Since by our work of Section 4 there will be no point of accumulation of epochs of jump times for the process, following Feller [(1966), Section X.3], we see [by first conditioning on $\{p(t)\}_{t\geq 0}$] that $\{F_{PP}(t,p), F_{PQ}(t,p), F_{QQ}(t,p)\}$ can be characterized as the unique minimal solution to the following system of differential equations (we use $\dot{F}$ to denote the derivative of $F$ with respect to $t$):

$$\begin{aligned}
\dot{F}_{PP} &= \frac{1-F_{PP}}{2p} + \left(\frac{\mu_2 q}{p} + rq\right)(F_{PQ} - F_{PP}) \\
&\quad + \frac{1}{2}(-\mu_1 p + \mu_2 q + spq)F'_{PP} + \frac{1}{4}pqF''_{PP}, \\
\dot{F}_{PQ} &= \frac{1}{2}\left(\frac{p\mu_1}{q} + rp\right)(F_{PP} - F_{PQ}) + \frac{1}{2}\left(\frac{q\mu_2}{p} + rq\right)(F_{QQ} - F_{PQ}) \\
&\quad + \frac{1}{2}(-\mu_1 p + \mu_2 q + spq)F'_{PQ} + \frac{1}{4}pqF''_{PQ}, \\
\dot{F}_{QQ} &= \frac{1-F_{QQ}}{2q} + \left(\frac{\mu_1 p}{q} + rp\right)(F_{PQ} - F_{QQ}) \\
&\quad + \frac{1}{2}(-\mu_1 p + \mu_2 q + spq)F'_{QQ} + \frac{1}{4}pqF''_{QQ}.
\end{aligned} \tag{23}$$

Suppose that the mutation rate to a novel allele at the neutral site is $\nu$, corresponding to the rescaling $\nu \mapsto \nu/N$ of the model of Section 2, then conditional on the individuals having coalesced at time $t$, the probability that they are identical in state, that is, that there has been no mutation since time $t$ along either of their lines of descent, is $e^{-2\nu t}$. Conditioning on the time to coalescence then gives

$$f_{PP}(p) = \int_0^\infty e^{-2\nu t}\frac{dF_{PP}(t,p)}{dt}\,dt$$

with similar expressions for $f_{PQ}(p)$ and $f_{QQ}(p)$. Integration by parts then shows that, under the diffusion approximation, the probabilities of identity satisfy

$$0 = -2\nu f_{PP} + \frac{1-f_{PP}}{2p} + \left(\frac{\mu_2 q}{p} + rq\right)(f_{PQ} - f_{PP})$$

$$
\begin{aligned}
&+\frac{1}{2}(-\mu_1 p+\mu_2 q+spq)f'_{PP}+\frac{1}{4}pqf''_{PP},\\
0&=-2\nu f_{PQ}+\frac{1}{2}\left(\frac{p\mu_1}{q}+rp\right)(f_{PP}-f_{PQ})\\
&+\frac{1}{2}\left(\frac{q\mu_2}{p}+rq\right)(f_{QQ}-f_{PQ})\\
&+\frac{1}{2}(-\mu_1 p+\mu_2 q+spq)f'_{PQ}+\frac{1}{4}pqf''_{PQ},\\
0&=-2\nu f_{QQ}+\frac{1-f_{QQ}}{2q}+\left(\frac{\mu_1 p}{q}+rp\right)(f_{PQ}-f_{QQ})\\
&+\frac{1}{2}(-\mu_1 p+\mu_2 q+spq)f'_{QQ}+\frac{1}{4}pqf''_{QQ}.
\end{aligned}
\tag{24}
$$

We now identify the probabilities of identity in state as the *minimal* solution to this system. Again following Feller [(1969), Section X.3], $\{F_{PP}(t,p)$, $F_{PQ}(t,p)$, $F_{QQ}(t,p)\}$, the minimal solution to (23) is most easily constructed via an iterative procedure. At the $n$th stage of the iteration, the functions $\{F^{(n)}_{PP}(t,p)$, $F^{(n)}_{PQ}(t,p)$, $F^{(n)}_{QQ}(t,p)\}$ are obtained by conditioning the number of jumps that the process can make by time $t$ to be at most $n$. Since (by our work of Section 4) the total number of jumps that the process can make by time $t$ is finite, as $n\to\infty$ this sequence of functions really does converge to the distribution function of the coalescence times. If we now define

$$f^{(n)}_{PP}=\int_0^\infty e^{-2\nu t}\frac{dF^{(n)}(t,p)}{dt}\,dt,$$

with parallel definitions for $f^{(n)}_{PQ}$ and $f^{(n)}_{QQ}$, then as $n\to\infty$, $\{f^{(n)}_{PP}(p)$, $f^{(n)}_{PQ}(p)$, $f^{(n)}_{QQ}(p)\}$ converges to the minimal solution to the system of equations (24).

Combining the above yields the following.

THEOREM 6.1. *Under the diffusion approximation, if the process of allele frequencies is assumed to have reached stationarity, then the probabilities of identity in state for a sample of size two whose types at the selected site are known, denoted* $\{f_{PP}(p), f_{PQ}(p), f_{QQ}(p)\}$, *are given by the minimal solution to the system of equations* (24).

Before exploring the equations numerically, we make precise the iteration that we used above to construct the minimal solution. The sequence of functions $\{f^{(n)}_{PP}, f^{(n)}_{PQ}, f^{(n)}_{QQ}\}_{n\geq 0}$ is obtained as follows. First set $(f^{(0)}_{PP}(p)$, $f^{(0)}_{PQ}(p)$,

$f_{QQ}^{(0)}(p)) \equiv (0,0,0)$ for $p \in [0,1]$. Then for $n \geq 1$,

$$\begin{aligned}
0 &= -2\nu f_{PP}^{(n)} + \frac{1-f_{PP}^{(n)}}{2p} + \left(\frac{\mu_2 q}{p} + rq\right)(f_{PQ}^{(n-1)} - f_{PP}^{(n)}) \\
&\quad + \frac{1}{2}(-\mu_1 p + \mu_2 q + spq)\frac{d}{dp}f_{PP}^{(n)} + \frac{1}{4}pq\frac{d^2}{dp^2}f_{PP}^{(n)}, \\
0 &= -2\nu f_{PQ}^{(n)} + \frac{1}{2}\left(\frac{p\mu_1}{q} + rp\right)(f_{PP}^{(n-1)} - f_{PQ}^{(n)}) \\
&\quad + \frac{1}{2}\left(\frac{q\mu_2}{p} + rq\right)(f_{QQ}^{(n-1)} - f_{PQ}^{(n)}) \\
&\quad + \frac{1}{2}(-\mu_1 p + \mu_2 q + spq)\frac{d}{dp}f_{PQ}^{(n)} + \frac{1}{4}pq\frac{d^2}{dp^2}f_{PQ}^{(n)}, \\
0 &= -2\nu f_{QQ}^{(n)} + \frac{1-f_{QQ}^{(n)}}{2q} + \left(\frac{\mu_1 p}{q} + rp\right)(f_{PQ}^{(n-1)} - f_{QQ}^{(n)}) \\
&\quad + \frac{1}{2}(-\mu_1 p + \mu_2 q + spq)\frac{d}{dp}f_{QQ}^{(n)} + \frac{1}{4}pq\frac{d^2}{dp^2}f_{QQ}^{(n)}.
\end{aligned}$$

As $n \to \infty$, the functions $f_{PP}^{(n)}$, $f_{PQ}^{(n)}$ and $f_{QQ}^{(n)}$ converge (monotonically) to the minimal solution of (24).

Since the system (24), arose by integrating the Kolmogorov backward equations (23), the boundary conditions are implicitly prescribed. However, in order to solve the equations numerically, we require explicit expressions. The first thing that we must check is that

$$p\frac{df_{PP}^{(n)}}{dp}, \quad p\frac{df_{PQ}^{(n)}}{dp}, \quad p\frac{d^2 f_{PP}^{(n)}}{dp^2}, \quad p\frac{d^2 f_{PQ}^{(n)}}{dp^2}$$

all tend to 0 as $p$ tends to 0 and similarly,

$$(1-p)\frac{df_{PQ}^{(n)}}{dp}, \quad (1-p)\frac{df_{QQ}^{(n)}}{dp}, \quad (1-p)\frac{d^2 f_{PQ}^{(n)}}{dp^2}, \quad (1-p)\frac{d^2 f_{QQ}^{(n)}}{dp^2}$$

all tend to 0 as $p$ tends to 1. The method is lengthy, but completely standard. For each equation, first use the Frobenius method of solution in series to find two linearly independent solutions to the corresponding homogeneous equation. (One solution will be singular at $p=0$ and the other will not.) Then use the method of variation of parameters to write down the corresponding Green's function and finally integrate to obtain the solution to the original (inhomogeneous) equation. We omit the details. They can be found in any standard text on ordinary differential equations, for example, Simmons (1974).

Granted the above, we can now read off the boundary conditions from the equations. Letting $p$ tend to 0 in the first equation and to 1 in the last equation of (24) gives

$$f_{PP}^{(n)}(0) = \frac{1 + 2\mu_2 f_{PQ}^{(n-1)}(0)}{1 + 2\mu_2}, \qquad f_{QQ}^{(n)}(1) = \frac{1 + 2\mu_1 f_{PQ}^{(n-1)}(1)}{1 + 2\mu_1}.$$

Letting $p$ tend to 1 in the first equation and to 0 in the last equation of (24) we see that we must have

$$(1 + 4\nu) f_{PP}^{(n)}(1) + \mu_1 \frac{df_{PP}^{(n)}}{dp}(1) = 1, \qquad (1 + 4\nu) f_{QQ}^{(n)}(0) - \mu_2 \frac{df_{QQ}^{(n)}}{dp}(0) = 1.$$

The second equation yields

$$f_{PQ}^{(n)}(0) = f_{PP}^{(n-1)}(0), \qquad f_{PQ}^{(n)}(1) = f_{QQ}^{(n-1)}(1).$$

REMARK 6.2. Our system of differential equations is for a sample of size two from a population that can be in just two possible states. Clearly this is a very special situation. Of course it is readily extended to larger systems. However, to characterize the transition probabilities for a sample of size $N$ from a population with $m$ possible states requires $\sum_{n=1}^{N} \sum_{k=1}^{m} \binom{m}{k} \binom{n-1}{k-1}$ equations. The distribution of the time to the most recent common ancestor of the sample requires $\sum_{n=2}^{N} \sum_{k=1}^{m} \binom{m}{k} \binom{n-1}{k-1}$ equations. The probability of identity in state for a sample of size two from a population distributed amongst $m$ genetic backgrounds requires $m(m+3)/2$ equations. Evidently for large samples or complex genetic backgrounds the approach will become intractable.

**7. Numerical examples.** In this section we illustrate the accuracy of (24) when compared to a direct solution of the matrix equations (1) and (2) that gave us the exact probabilities of identity for the Wright–Fisher model. We then use the equations to illustrate the potentially important influence of the fluctuations on the probabilities of identity in allelic state at the neutral locus. We concentrate exclusively on the case of balancing selection $s = s_0(p_0 - p)$ with $p_0 = 1/2$. We also set the recombination rate $r = 0$.

These results could all be obtained for low selection rates using the ancestral selection graph methods in which lineages branch into three potential ancestors at rate $s_0$. Here we used Mathematica to solve the differential equations numerically (a process that takes only seconds of computer time irrespective of the strength of selection). In all the figures we have included a plot of the stationary distribution for the allele frequencies for the parameter values used. The unique stationary distribution for the process has density

$$m(p)=\beta p^{2\mu_2-1}(1-p)^{2\mu_1-1}\exp(-\tfrac{1}{2}s_0(p^2+(1-p)^2)),\qquad p\in(0,1),$$

where the constant $\beta$ is chosen so that $\int_0^1 m(p)\,dp=1$. Note that when the selection is very strong, rounding errors mean that the explosion of the density at the margins is not visible on the plot.

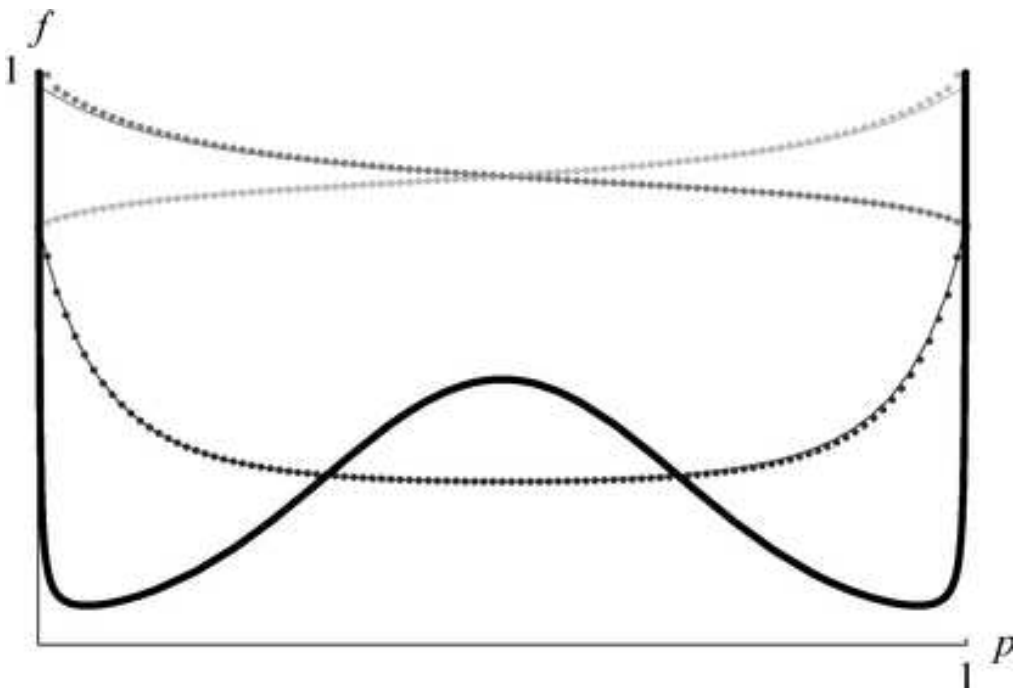

Fig. 1. *Comparison of solutions to the matrix equations for identity in state obtained from the Wright–Fisher model to those of the differential equations obtained from the diffusion approximation. The three thin lines are (in descending order at $p=0$) $f_{PP}$, $f_{QQ}$ and $f_{PQ}$. The dotted lines are the solutions to the matrix equations and the bold line is the stationary distribution of $p$. The parameter values are $N=50$, $s_0=0.16$, $p_0=0.5$, $\mu_1=0.0005=\mu_2$, $\nu=0.002$.*

Figure 1 compares the solution to the matrix equations (1) and (2) to the diffusion approximation obtained by solving the system (24). Even for the modest population size (fifty diploid individuals), the accuracy of the diffusion approximation is striking. We obtained similar results with other parameter values. In this example the mutation rate between selected loci was chosen to be very small as this is the case when one expects the diffusion approximation to be most likely to break down.

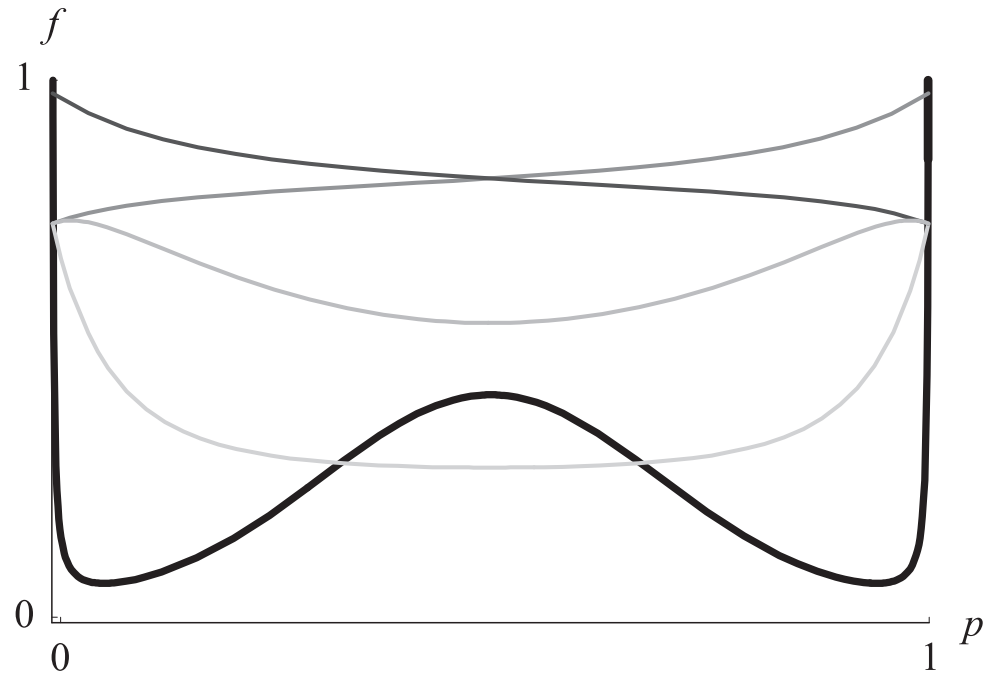

Fig. 2. *From top to bottom at $p=0$ the plotted functions correspond to $f_{PP}$, $f_{QQ}$, $\overline{f}$ and $f_{PQ}$. The thick line is the stationary distribution. The parameter values are $\nu=0.1$, $\mu_1=0.025=\mu_2$, $p_0=0.5$, $s_0=0.16$.*

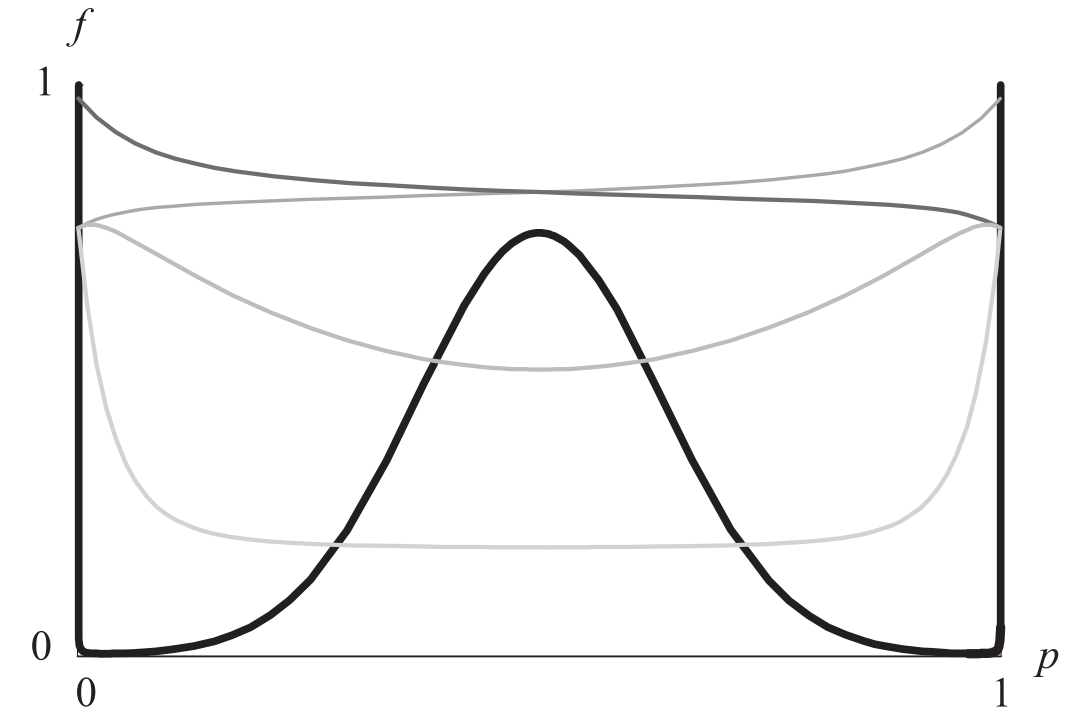

FIG. 3. *The same as Figure 2 except that now $s_0 = 0.32$.*

In Figures 2 and 3 the probabilities of identity are plotted for the case of strong balancing selection. In both cases the mean value of the allele frequencies at the selected locus when the population has reached stationarity is 1/2. The probability $\overline{f} = p^2 f_{PP} + 2pqf_{PQ} + q^2 f_{QQ}$ of identity for two individuals selected at random is also plotted. This function should be compared to the constant value 0.43 that one obtains by setting $p \equiv \frac{1}{2}$ and using the standard structured coalescent model. Even when the strength of selection is rather strong, this is a poor approximation. In Figure 4 we plot the result of integrating the function $\overline{f}$ against the stationary distribution of the allele frequencies for different strengths of selection (all other parameters being as in Figures 2 and 3). As we see the strength of selection has to be very strong indeed before the value predicted by the standard structured coalescent can be regarded as a good approximation. Finally, in Figure 5 we compare the mean time to coalescence for a sample of size two as a function of the strength of selection.

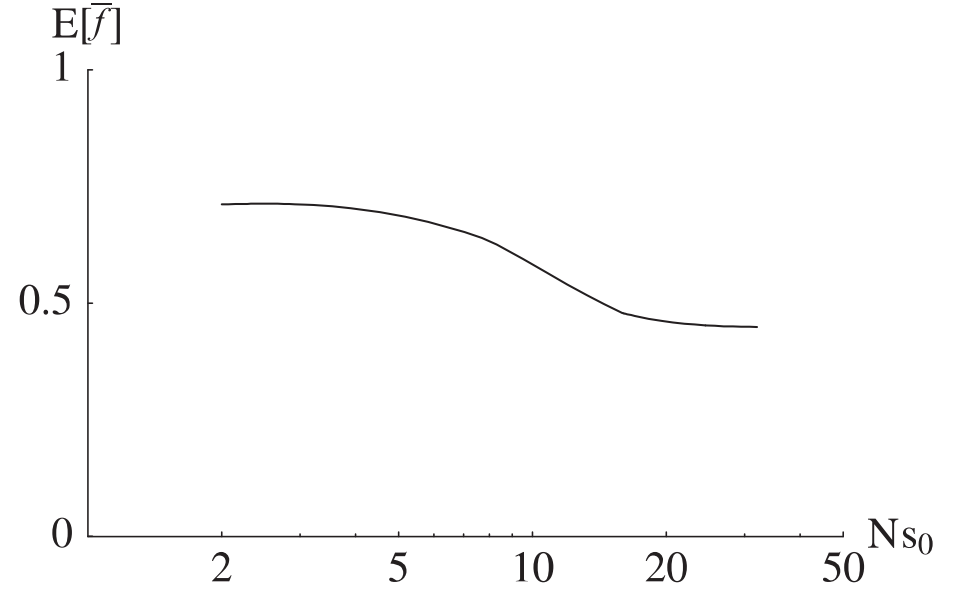

FIG. 4. *Prediction of the diffusion approximation for the integral of $\overline{f}$ against the stationary distribution for the allele frequencies at the selected locus as a function of the strength of selection. All other parameters are chosen as in Figures 2 and 3.*

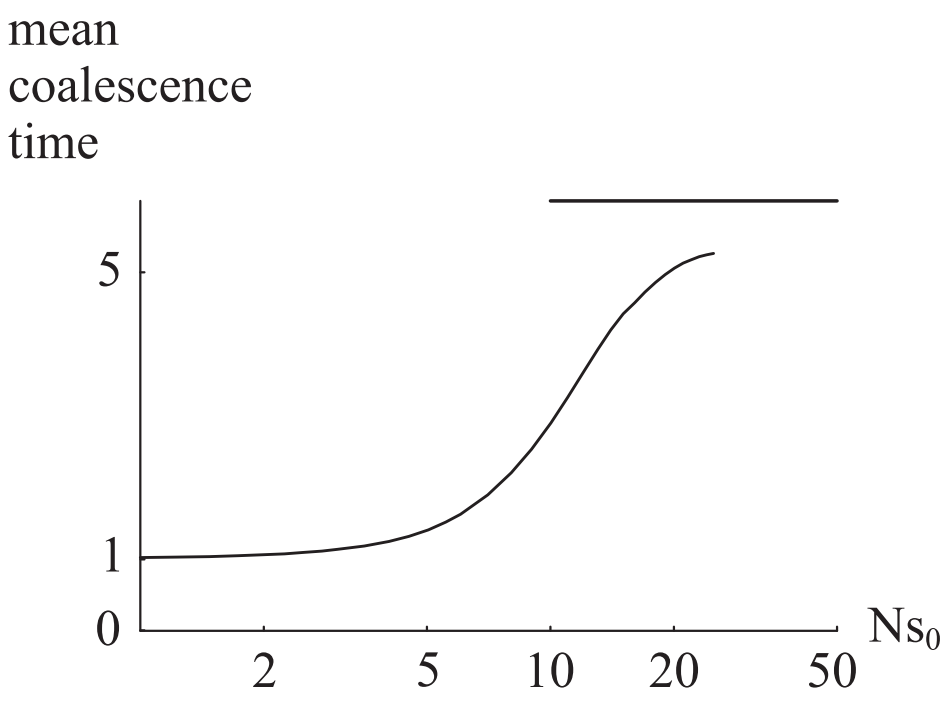

FIG. 5. *Predictions of the diffusion approximation for mean coalescence time scaled relative to the neutral expectation (2N) for a sample of size two. The parameters are as in Figure 4. The thick line on the upper right is the prediction from the structured coalescent (6N).*

We refer to the companion paper, Barton and Etheridge (2004), for a more detailed investigation and discussion of the biological issues raised.

**Acknowledgments.** Much of this work was completed when the second author was visiting the University of Edinburgh. She would like to thank everyone there for their hospitality. The third author would like to thank M. Nordborg for bringing the problem to her attention and R. C. Griffiths for helpful discussions.

N. H. Barton
Institute of Cell, Animal
and Population Biology
University of Edinburgh
King's Buildings
West Mains Road
Edinburgh EH9 3JT
United Kingdom
e-mail: n.barton@ed.ac.uk

A. M. Etheridge
Department of Statistics
University of Oxford
1 South Parks Road
Oxford OX1 3TG
United Kingdom
e-mail: etheridg@stats.ox.ac.uk

A. K. Sturm
Institute of Mathematics
Technical University Berlin
Strasse des 17. Juni 136
D-10623 Berlin
Germany
e-mail: sturm@math.tu-berlin.de